\documentclass[11pt]{amsart}

\swapnumbers

\usepackage{amssymb,amsfonts}
\usepackage{euscript}

\usepackage[ps,matrix,arrow,curve,cmtip]{xy}

\title{A streamlined proof of Goodwillie's $n$-excisive approximation}
\author{Charles Rezk}
\date{ \today}
\address{Department of Mathematics \\
University of Illinois at Urbana-Champaign \\ 
Urbana, IL}
\email{rezk@math.uiuc.edu}
\urladdr{http://www.math.uiuc.edu/~rezk}

\numberwithin{equation}{section}

\makeatletter
  \let\c@subsection\c@equation
\makeatother

\theoremstyle{plain}   

\newtheorem{lemma}[subsection]{Lemma}

\theoremstyle{remark}

\theoremstyle{plain}

\newcommand{\ra}{\rightarrow}

\newcommand{\xra}{\xrightarrow}
\newcommand{\la}{\leftarrow}

\DeclareMathOperator{\hocolim}{hocolim}
\DeclareMathOperator{\holim}{holim}

\begin{document}

\begin{abstract}
We give a shorter proof of Lemma 1.9 from Goodwillie, ``Calculus
III'', which is the key step in proving that the construction $P_nF$
gives an $n$-excisive functor.
\end{abstract}

\maketitle


\newcommand{\spaces}{\mathcal{U}}
\newcommand{\poset}{\mathcal{P}}
\newcommand{\cubeX}{\mathcal{X}}

\section{Introduction}

For a homotopy functor $F$ from spaces to spaces, Goodwillie has
defined the notion of an ``$n$-excisive approximation'', which is a
homotopy functor $P_nF$ together with a natural transformation
$p_nF\colon F\ra P_nF$.  In \cite[Thm.\ 1.8]{goodwillie-calculus-3} it
is shown that the functor $P_nF$ is in fact an $n$-excisive functor,
and therefore that $p_nF$  is the universal example of a map from
$F$ to an $n$-excisive functor.  

Earlier work of
Goodwillie had shown that $P_nF$ is $n$-excisive under additional
hypotheses involving connectivity.  
The notable feature of the proof given in \cite{goodwillie-calculus-3}
is that that no hypotheses involving connectivity are needed.  In
fact, the argument is entirely general, and will work in any homotopy
theory in which directed homotopy colimits commute with finite
homotopy limits.

Goodwillie's
proof relies the following  ``lemma'' \cite[Lemma.\
1.9]{goodwillie-calculus-3}.  (The notions of ``cartesian'' and
``strongly cocartesian cube'' are defined in \S1 of
\cite{goodwillie-calculus-3}.  The definitions of $T_nF$ and $t_nF$
are given below.)  
\begin{lemma}\label{lemma}
Let $\cubeX$ be any strongly cocartesian $n$-cube in $\spaces$, and
let $F$ be any homotopy functor.  The map of cubes
$(t_{n}F)(\cubeX)\colon F(\cubeX)\ra (T_{n}F)(\cubeX)$ factors
through some cartesian cube.
\end{lemma}
Goodwillie's 
proof of this lemma 
is, as he notes, ``a little 
opaque''.  In fact, though the proof gives an explicit factorization of
$(t_nF)(\cubeX)$ through a cartesian cube, the cube in question is
difficult to describe, and does not seem to play any natural role.

The purpose of this note is to give a much simplified proof of
Goodwillie's lemma (though in the same spirit as Goodwillie's), and
thus a simplified proof of the construction of the $n$-excisive
approximation.  We will assume that the reader is familiar with
\cite{goodwillie-calculus-3}, and we assume the context and notation
of \S1 of that paper.

The author was supported under NSF grants DMS-0505056 and DMS-1006054.

\section{Proof of the lemma}

Let $\poset(n)$
denote the poset of subsets of $\{1,\cdots,n\}$, and let
$\poset_0(n)\subset \poset(n)$ be the poset of non-empty subsets.

If $F\colon \mathcal{C}\ra \mathcal{D}$ is a homotopy functor, Goodwillie
defines a functor 
$T_{n}F\colon \mathcal{C}\ra \mathcal{D}$ and natural map $t_nF\colon F\ra
T_{n-1}F$ by
\[
F(X) \xra{t_{n}F} \holim_{U\in \poset_0(n+1)} F(X*U).
\]

\begin{proof}[Proof of \ref{lemma}]
We write $n$ instead of $n+1$.
Given any cube $\cubeX$  and a set $U\in \poset(n)$, define a cube
$\cubeX_U$ by
\[
\cubeX_U (T) = \hocolim \left( \cubeX(T) \la \coprod_{s\in U} \cubeX(T) \ra
\coprod_{s\in U} \cubeX(T\cup\{s\}) \right).
\]
We have $\cubeX_\varnothing(T)\approx \cubeX(T)$, and there is an evident map
$\alpha\colon \cubeX_U(T) \ra \cubeX(T)*U$, which 
is natural in both $T$ and $U$.  

The map $(t_{n-1}F)(\cubeX)$ factors as follows:
\[
F(\cubeX(T)) \ra \holim_{U\in \poset_0(n)} F(\cubeX_U(T)) \ra \holim_{U\in
  \poset_0(n)} F(\cubeX(T)*U) \approx (T_{n-1}F)(\cubeX(T)).
\]

Now suppose that $\cubeX$ is strongly cocartesian.  Then there are
natural weak equivalences $\cubeX_U(T) \approx \cubeX(T\cup
U)$.  The maps $\cubeX(T\cup U)\ra
\cubeX(T\cup\{s\} \cup U)$ are isomorphisms for $s\in U$, and
thus if $U$ is non-empty the cube $T\mapsto F(\cubeX_U(T))$ is
cartesian.  Therefore
$\holim_{U\in \poset_0(n)} F(\cubeX_U(T))$ is a homotopy limit of
cartesian cubes, and thus is cartesian.
\end{proof}

Note that this shows that if $T$ is non-empty, then
$U\mapsto F(\cubeX_U(T))$ is cartesian, so
that $F(\cubeX(T))\ra 
\holim_{U\in\poset_0(n)} F(\cubeX_U(T))$ is a weak equivalence for
$T\neq\varnothing$.  For $T=\varnothing$, we see that
$\holim_{U\in\poset_0(n)} F(\cubeX_U(\varnothing))\approx \holim_{U\in
  \poset_0(n)}F(\cubeX(U))$.


\begin{thebibliography}{Goo03}

\bibitem[Goo03]{goodwillie-calculus-3}
Thomas~G. Goodwillie, \emph{Calculus. {III}. {T}aylor series}, Geom. Topol.
  \textbf{7} (2003), 645--711 (electronic). \MR{MR2026544}

\end{thebibliography}
\newcommand{\noopsort}[1]{} \newcommand{\printfirst}[2]{#1}
  \newcommand{\singleletter}[1]{#1} \newcommand{\switchargs}[2]{#2#1}
\providecommand{\bysame}{\leavevmode\hbox to3em{\hrulefill}\thinspace}
\providecommand{\MR}{\relax\ifhmode\unskip\space\fi MR }
\providecommand{\MRhref}[2]{%
  \href{http://www.ams.org/mathscinet-getitem?mr=#1}{#2}
}
\providecommand{\href}[2]{#2}

\end{document}